\documentclass[preprint,12pt]{elsarticle}

\usepackage[T1]{fontenc}
\usepackage[utf8]{inputenc}
\usepackage{amsmath,amssymb,amsthm}
\usepackage{hyperref}
\usepackage{tikz}

\newtheorem{theorem}{Theorem}[section]
\newtheorem{lemma}[theorem]{Lemma}
\newtheorem{proposition}[theorem]{Proposition}
\newtheorem{corollary}[theorem]{Corollary}
\newtheorem{definition}[theorem]{Definition}
\newtheorem{remark}[theorem]{Remark}

\newcommand{\CVS}{\mathrm{CVS}}
\newcommand{\Tl}{Tl}
\newcommand{\R}{\mathbb{R}}

\journal{Journal Name}  % 投稿时填入目标期刊

\begin{document}

\begin{frontmatter}

\title{Dimensionally absolutely stable but non-diagonalizable \(T\)-meshes}

\author[snnu]{Mao Shi\corref{cor1}}
\ead{shimao@snnu.edu.cn}  % 请替换为实际邮箱
\affiliation[snnu]{organization={School of Mathematics and Statistics, Shaanxi Normal University},
            addressline={},
            city={Xi'an},
            postcode={},
            state={Shaanxi},
            country={China}}

\cortext[cor1]{Corresponding author}

\begin{abstract}
We construct an infinite family of rectangular \(T\)-meshes whose conformality-vector-space dimension is constant over all structurally similar geometric realizations, although the meshes are not diagonalizable. For every integer \(d\ge 2\), the construction has \(2(d+2)\) \(\Tl\)-edges and intersection graph \(K_{d+2,d+2}\). Every \(\Tl\)-edge contains at least \(d+2\) vertices, so no edge is vanishable under the usual vertex-count criterion for bidegree \((d,d)\) and smoothness \((d-1,d-1)\). We prove that the conformality-vector space is one-dimensional throughout the structural similarity class. Its generator is the outer product of two vectors of barycentric interpolation weights, extended by four zero coordinates. Non-diagonalizability follows by examining the last edge in any proposed diagonalizing ordering. Thus, for the geometric structural similarity class specified here, diagonalizability is not a necessary condition for dimensional absolute stability.
\end{abstract}


\end{frontmatter}

\section{Introduction}

The dimension of a spline space over a \(T\)-mesh is closely related to the rank of its conformality matrix. For splines of bidegree \((d,d)\) and smoothness \((d-1,d-1)\), the constraints associated with a \(\Tl\)-edge can be written as polynomial moment equations. Their coefficients depend on the coordinates of the vertices on that edge. This is part of a broader literature on dimension instability and stability of spline spaces over \(T\)-meshes; see Li and Chen \cite{liChen2011}, Berdinsky et al. \cite{berdinsky2012}, and Guo et al. \cite{guo2015}. The smoothing cofactor-conformality method is developed in Li and Deng \cite{liDeng2016}. Stability properties for highest-order smoothness spaces over \(T\)-meshes, including the notion of dimensional absolute stability, have been studied by Huang and Chen \cite{huang2024,huang2026}.

Diagonalizability is an ordering property of the \(\Tl\)-edges. It provides a combinatorial mechanism for controlling these constraints. The converse question is whether coordinate-independent dimension, even throughout a structural similarity class, forces such an ordering to exist.

This paper gives a negative answer for the geometric structural similarity class defined below. For each \(d\ge 2\), we construct a rectangular \(T\)-mesh with no vanishable \(\Tl\)-edges whose conformality-vector-space dimension is identically one throughout that class, but for which no diagonalizing ordering exists.

The construction uses a complete bipartite intersection pattern. This pattern forces a shared set of coordinates on all edges of each orientation. The resulting compatibility equations admit a rank-one factorization involving classical barycentric interpolation weights.

To make the scope of the result precise, we explicitly state the structural similarity and stability conventions used here. In particular, the realizations under consideration are geometric collections of horizontal and vertical segments in the plane, not abstract incidence systems with independently assigned local coordinates.

\section{Definitions and conventions}

Fix an integer \(d\ge 2\). The associated spline setting is bidegree \((d,d)\) with smoothness \((d-1,d-1)\).

A rectangular \(T\)-mesh is a partition of a rectangle into axis-aligned rectangles, allowing \(T\)-junctions. An interior \(l\)-edge is a maximal connected straight union of interior mesh edges. A cross-cut has both endpoints on the boundary of the domain. A \(\Tl\)-edge has both endpoints at interior \(T\)-junctions.

For an \(l\)-edge \(\ell\), let \(\operatorname{Vert}(\ell)\) be its set of vertices, including its endpoints, and write
\[
n(\ell)=\#\operatorname{Vert}(\ell).
\]
A \(\Tl\)-edge is called vanishable if
\[
n(\ell)\le d+1.
\]
A \(T\)-connected component is a connected collection of \(\Tl\)-edges, where connectivity is determined by shared vertices. Within such a component, a vertex belonging to two component edges is called a multi-vertex, and a vertex belonging to exactly one component edge is called a mono-vertex.

\subsection{Geometric structural similarity}

The following two definitions adapt the notion of a generalized \(T\)-connected component and the structurally similar map of Huang and Chen \cite{huang2026} to the planar geometric setting used in this paper.

\begin{definition}[Geometric realization]
A geometric component consists of finitely many horizontal and vertical closed line segments in \(\R^2\), together with a finite set of distinct marked vertices on each segment. Endpoints and intersections of component segments are marked.

Segments of the same orientation do not share a vertex. Segments of opposite orientations meet in at most one vertex. Every marked vertex belongs to one or two component segments. The segment intersection graph is connected.
\end{definition}

The segments in Definition~2.1 represent the \(l\)-edges of a generalized geometric component. Additional marked vertices represent mono-vertices. The common-coordinate requirement in Definition~2.1 (same orientation segments do not share a vertex) is essential for the analysis below; it is precisely the planar geometric counterpart of the common-coordinate property discussed by Huang and Chen \cite{huang2026}.

\begin{definition}[Structural similarity]
Two geometric components \(\mathcal C\) and \(\widetilde{\mathcal C}\) are structurally similar if there is a bijection
\[
\phi:\mathcal E(\mathcal C)\longrightarrow \mathcal E(\widetilde{\mathcal C})
\]
between their edge sets such that
\[
n(\phi(\ell))=n(\ell)
\]
for every edge \(\ell\), and
\[
\ell\cap\ell'\ne\emptyset
\iff
\phi(\ell)\cap\phi(\ell')\ne\emptyset
\]
for every pair of distinct edges.
\end{definition}

This is the geometric counterpart of the structurally similar map of Huang and Chen \cite{huang2026}: the condition \(\ell\cap\ell'\ne\emptyset\iff\phi(\ell)\cap\phi(\ell')\ne\emptyset\) plays the role of the incidence-preservation condition \(\beta(\ell_i\cap\ell_j)=\beta(\ell_i)\cap\beta(\ell_j)\) in the original definition.

The class of all such geometric realizations is denoted by \([\mathcal C]_{\mathrm S}\).

Structural similarity does not prescribe numerical coordinates, distances, or the order of mono-vertices. It also permits an exchange of the two coordinate directions.

\subsection{Conformality-vector spaces}

Assign one real scalar \(\gamma_v\) to each distinct vertex \(v\) of a component \(\mathcal C\). On a horizontal edge, let \(z_v\) be the horizontal coordinate of \(v\); on a vertical edge, let \(z_v\) be its vertical coordinate.

The moment conditions on an edge \(\ell\) are
\begin{equation}
\sum_{v\in\operatorname{Vert}(\ell)} z_v^q\gamma_v=0,
\qquad q=0,\ldots,d.
\label{eq:moment}
\end{equation}
At a multi-vertex, the same variable is used in the conditions for both incident component edges.

The conformality-vector space is
\begin{equation}
\CVS_d(\mathcal C)
=
\left\{
(\gamma_v)_{v\in\operatorname{Vert}(\mathcal C)}
:\ \eqref{eq:moment}\text{ holds on every edge}
\right\}.
\label{eq:cvs}
\end{equation}
Equivalently, \(\CVS_d(\mathcal C)=\ker M_d(\mathcal C)\), where \(M_d(\mathcal C)\) is the matrix of these moment conditions.

\begin{definition}[Dimensional absolute stability; Huang and Chen \cite{huang2026}]
A component \(\mathcal C\) is dimensionally absolutely stable if
\(\dim\CVS_d(\widetilde{\mathcal C})\) is independent of
\(\widetilde{\mathcal C}\in[\mathcal C]_{\mathrm S}\).
\end{definition}

A \(T\)-mesh with a single \(T\)-connected component is called dimensionally absolutely stable if that component has this property.

All meshes constructed below have exactly one \(T\)-connected component.

\subsection{Diagonalizability}

\begin{definition}[Diagonalizability; Li and Deng \cite{liDeng2016}, see also Huang and Chen \cite{huang2026}, Definition~2.3]
A component with edges \(\ell_1,\ldots,\ell_t\) is diagonalizable if its edges admit an ordering such that
\begin{equation}
\#
\left(
\operatorname{Vert}(\ell_i)
\setminus
\bigcup_{j<i}\operatorname{Vert}(\ell_j)
\right)
\ge d+1,
\qquad i=1,\ldots,t.
\label{eq:diag}
\end{equation}
\end{definition}

A \(T\)-mesh with one \(T\)-connected component is diagonalizable if its component is diagonalizable.

\begin{remark}[Scope of the definitions]
The common-coordinate property of a planar geometric realization is essential below. The argument does not assert stability for a different model in which each edge is assigned its own independent set of local vertex coordinates.

Accordingly, applying the result to an externally formulated conjecture requires matching that conjecture's definitions with Definitions~2.1--2.4.
\end{remark}

\section{A classical Vandermonde identity}

We use the following standard consequence of the divided difference formula and the Lagrange interpolation formula; see de Boor \cite{deboor2005} and Berrut and Trefethen \cite{berrut2004}. The explicit vector consists of the barycentric interpolation weights. We state the result without proof.

\begin{lemma}[One-dimensional Vandermonde null space]
Let \(r\ge 2\), and let \(z_1,\ldots,z_r\) be pairwise distinct real numbers. Set
\begin{equation}
V(z_1,\ldots,z_r)
=
\begin{pmatrix}
1&\cdots&1\\
z_1&\cdots&z_r\\
\vdots&\ddots&\vdots\\
z_1^{r-2}&\cdots&z_r^{r-2}
\end{pmatrix}
\in\R^{(r-1)\times r}.
\label{eq:V}
\end{equation}
Then \(\operatorname{rank}V=r-1\) and \(\ker V\) is one-dimensional.

More explicitly, if
\[
P(z)=\prod_{k=1}^r(z-z_k),
\]
then
\[
w_j=\frac{1}{P'(z_j)},\qquad j=1,\ldots,r,
\]
defines a generator \(w=(w_1,\ldots,w_r)^T\) of \(\ker V\). Every component of \(w\) is nonzero. Consequently, every nonzero vector in \(\ker V\) has no zero component.
\end{lemma}

\section{Construction and non-diagonalizability}

Let
\[
r=d+2
\]
and consider the rectangle
\[
\Omega=[0,r+3]\times[0,r+1].
\]
Insert the two vertical cross-cuts
\[
C_{\mathrm L}=\{1\}\times[0,r+1],
\qquad
C_{\mathrm R}=\{r+2\}\times[0,r+1].
\]
Insert the \(r\) vertical segments
\begin{equation}
V_j=\{j+1\}\times[1,r],
\qquad j=1,\ldots,r,
\label{eq:Vj}
\end{equation}
the two horizontal segments
\begin{equation}
H_1=[1,r+2]\times\{1\},
\qquad
H_r=[1,r+2]\times\{r\},
\label{eq:H1Hr}
\end{equation}
and the remaining horizontal segments
\begin{equation}
H_i=[2,r+1]\times\{i\},
\qquad i=2,\ldots,r-1.
\label{eq:Hi}
\end{equation}
Together with the boundary of \(\Omega\), these segments define a \(T\)-mesh denoted by \(\mathcal T_d\).

\begin{proposition}
The construction defines a rectangular \(T\)-mesh. Its \(\Tl\)-edges are precisely
\[
H_1,\ldots,H_r,V_1,\ldots,V_r.
\]
They form one \(T\)-connected component \(\mathcal C_d\), whose intersection graph is \(K_{r,r}\). Moreover,
\begin{equation}
\begin{aligned}
n(H_1)&=n(H_r)=r+2,\\
n(H_i)&=r, && 2\le i\le r-1,\\
n(V_j)&=r, && 1\le j\le r.
\end{aligned}
\label{eq:counts}
\end{equation}
In particular, no \(\Tl\)-edge is vanishable.
\end{proposition}

\begin{proof}
The two cross-cuts partition \(\Omega\) into three vertical strips. The two outer strips are rectangles. Within the middle strip, \(H_1\) and \(H_r\) separate a lower rectangle, an upper rectangle, and a central rectangle.

The central rectangle is partitioned by the \(V_j\) into two side rectangles and \(r-1\) interior columns. The segments \(H_2,\ldots,H_{r-1}\) subdivide those interior columns into rectangles. Thus all mesh cells are rectangles.

Each \(H_1\) and \(H_r\) ends at \(T\)-junctions on the cross-cuts. Each \(V_j\) ends at \(T\)-junctions on \(H_1\) and \(H_r\). Each remaining \(H_i\) ends at \(T\)-junctions on \(V_1\) and \(V_r\). The listed segments are maximal in their respective directions.

Every \(H_i\) meets every \(V_j\), and no two listed edges of the same orientation meet. Hence the intersection graph is \(K_{r,r}\).

Every edge contains \(r\) multi-vertices. The two edges \(H_1\) and \(H_r\) additionally contain two mono-vertices each, at the cross-cuts. There are no other mono-vertices. The asserted counts follow.

Finally, \(r=d+2\), so every \(\Tl\)-edge contains more than \(d+1\) vertices.
\end{proof}

\begin{proposition}
For every \(d\ge 2\), the \(T\)-mesh \(\mathcal T_d\) is not diagonalizable.
\end{proposition}

\begin{proof}
Suppose that a diagonalizing ordering exists and consider its last edge \(\ell\). All other component edges precede \(\ell\).

If \(\ell\) is a vertical edge or one of \(H_2,\ldots,H_{r-1}\), all its vertices belong to preceding edges. Thus it contributes no new vertex.

If \(\ell\) is \(H_1\) or \(H_r\), all its multi-vertices belong to preceding vertical edges. Only its two mono-vertices remain.

In either case the number of remaining vertices is at most two. Since \(d\ge 2\), this is less than \(d+1\), contradicting \eqref{eq:diag}.
\end{proof}

\section{Absolute stability}

The essential point is to treat every member of the structural similarity class, not merely coordinate perturbations of the displayed mesh.

\begin{lemma}[Geometry of the similarity class]
Let \(\widetilde{\mathcal C}\in[\mathcal C_d]_{\mathrm S}\). After relabeling edges and, if necessary, exchanging the coordinate directions, its edges can be denoted by
\[
H_1,\ldots,H_r,V_1,\ldots,V_r
\]
with the following properties:
\begin{enumerate}
\item every \(H_i\) is horizontal and every \(V_j\) is vertical;
\item every \(H_i\) meets every \(V_j\);
\item \(H_1\) and \(H_r\) each contain two mono-vertices;
\item all remaining vertices are the \(r^2\) intersections \(H_i\cap V_j\).
\end{enumerate}
Furthermore, the supporting lines of the \(V_j\) have pairwise distinct horizontal coordinates \(x_j\), and the supporting lines of the \(H_i\) have pairwise distinct vertical coordinates \(y_i\).
\end{lemma}

\begin{proof}
By Proposition~4.1, the edge-incidence graph of \(\mathcal C_d\) is \(K_{r,r}\). Structural similarity preserves the edge-incidence relation, so \(\widetilde{\mathcal C}\) has the same incidence graph \(K_{r,r}\). Since \(K_{r,r}\) is connected, its bipartition is unique up to exchanging the two parts. In an axis-aligned geometric realization, the horizontal and vertical orientation classes give such a bipartition; hence, after possibly exchanging the coordinate directions, the two partite sets are precisely \(\{H_1,\ldots,H_r\}\) (horizontal) and \(\{V_1,\ldots,V_r\}\) (vertical). This proves properties~1 and~2.

Structural similarity also preserves the number of vertices on every corresponding edge. By \eqref{eq:counts}, in \(\mathcal C_d\) the two edges \(H_1,H_r\) have \(r+2\) vertices while all other edges have \(r\) vertices. Hence the same vertex counts hold in \(\widetilde{\mathcal C}\). Every edge meets exactly \(r\) edges of the opposite orientation, and each such intersection is a distinct vertex, so every edge has exactly \(r\) multi-vertices. An edge with \(r\) vertices therefore has no mono-vertex, whereas an edge with \(r+2\) vertices has exactly two mono-vertices. This proves properties~3 and~4.

Finally, by Definition~2.1, two same-orientation segments do not share a vertex. If two vertical edges had the same horizontal coordinate, then their intersections with any horizontal edge would coincide, making them share a vertex, a contradiction. Hence the \(x_j\) are pairwise distinct. The same argument applies to the \(y_i\).
\end{proof}

The labels \(H_1\) and \(H_r\) in Lemma~5.1 identify the two special edges; they need not be the geometrically lowest and highest edges of an arbitrary realization.

\begin{theorem}[Constant conformality-vector-space dimension]
For every \(\widetilde{\mathcal C}\in[\mathcal C_d]_{\mathrm S}\),
\[
\dim\CVS_d(\widetilde{\mathcal C})=1.
\]
More explicitly, let
\[
P_x(z)=\prod_{j=1}^r(z-x_j),
\qquad
P_y(z)=\prod_{i=1}^r(z-y_i),
\]
and define
\[
a_i=\frac{1}{P_y'(y_i)},
\qquad
b_j=\frac{1}{P_x'(x_j)}.
\]
Then every conformality vector has intersection coordinates
\begin{equation}
\gamma_{ij}=\mu a_i b_j,
\qquad 1\le i,j\le r,
\label{eq:gamma}
\end{equation}
for some \(\mu\in\R\), and all four mono-vertex coordinates are zero.
\end{theorem}

\begin{proof}
Use the labeling from Lemma~5.1. Set
\[
X=V(x_1,\ldots,x_r),
\qquad
Y=V(y_1,\ldots,y_r).
\]
Since \(d=r-2\), these are exactly the moment matrices on edges having no mono-vertices.

By Lemma~3.1,
\[
\ker X=\operatorname{span}\{b\},
\qquad
\ker Y=\operatorname{span}\{a\},
\]
and every component of \(a\) and \(b\) is nonzero.

Let \(\gamma_{ij}\) denote the coordinate associated with \(H_i\cap V_j\). Each \(V_j\) has only its \(r\) intersection vertices. Its moment conditions therefore give
\[
(\gamma_{1j},\ldots,\gamma_{rj})^T\in\ker Y.
\]
Consequently, for some \(\lambda_j\in\R\),
\[
\gamma_{ij}=a_i\lambda_j.
\]
Writing \(\lambda=(\lambda_1,\ldots,\lambda_r)^T\), the intersection-coordinate matrix is \(a\lambda^T\).

There is an ordinary horizontal edge, for example \(H_2\). Its moment conditions yield
\[
0=X(a_2\lambda)=a_2X\lambda.
\]
Since \(a_2\ne0\), it follows that \(\lambda\in\ker X\). Hence \(\lambda=\mu b\) for some \(\mu\in\R\), which proves \eqref{eq:gamma}. All remaining ordinary horizontal equations are then automatically satisfied.

Consider either special horizontal edge \(H_s\), where \(s\in\{1,r\}\). Let its mono-vertices have horizontal coordinates \(p_s\) and \(q_s\), with corresponding conformality coordinates \(\eta_s\) and \(\theta_s\). These are distinct vertices on a horizontal line, so \(p_s\ne q_s\).

For every \(k=0,\ldots,d\), the contribution of its intersection vertices is
\[
\sum_{j=1}^r x_j^k\gamma_{sj}
=
\mu a_s\sum_{j=1}^r x_j^k b_j
=0.
\]
Thus the remaining conditions are
\[
p_s^k\eta_s+q_s^k\theta_s=0,
\qquad k=0,\ldots,d.
\]
Taking \(k=0,1\) gives
\[
\begin{pmatrix}
1&1\\
p_s&q_s
\end{pmatrix}
\begin{pmatrix}
\eta_s\\
\theta_s
\end{pmatrix}
=0.
\]
The matrix is nonsingular because \(p_s\ne q_s\), so \(\eta_s=\theta_s=0\).

We have shown that every conformality vector has the asserted form. Conversely, \eqref{eq:gamma}, extended by four zero mono-vertex coordinates, satisfies every moment condition because \(Xa=0\) is not required: the relevant identities are \(Xb=0\) and \(Ya=0\). In particular, the horizontal equations use \(Xb=0\) and the vertical equations use \(Ya=0\).

Taking \(\mu=1\) gives a nonzero vector, and the solution space is therefore exactly one-dimensional.
\end{proof}

\begin{remark}
With a fixed vectorization convention, the full space is
\[
\CVS_d(\widetilde{\mathcal C})
=
\operatorname{span}
\left\{
\left(\operatorname{vec}(ab^T),0,0,0,0\right)
\right\}.
\]
The four final entries correspond to the mono-vertices.
\end{remark}

\begin{theorem}[An infinite counterexample family]
For every \(d\ge 2\), the rectangular \(T\)-mesh \(\mathcal T_d\) has no vanishable \(\Tl\)-edges, is dimensionally absolutely stable in the sense of Definition~2.3, and is not diagonalizable.
\end{theorem}

\begin{proof}
Proposition~4.1 gives the vertex-count condition. Theorem~5.2 gives dimensional absolute stability, and Proposition~4.2 gives non-diagonalizability.
\end{proof}

\begin{corollary}[Rank of the conformality matrix]
For every \(\widetilde{\mathcal C}\in[\mathcal C_d]_{\mathrm S}\), the matrix \(M_d(\widetilde{\mathcal C})\) has \(2r(r-1)\) rows, \(r^2+4\) columns, and rank
\[
\operatorname{rank}M_d(\widetilde{\mathcal C})=r^2+3.
\]
Its row-rank deficiency is
\[
2r(r-1)-(r^2+3)=(r-3)(r+1).
\]
\end{corollary}

\begin{proof}
There are \(2r\) edges, each contributing \(d+1=r-1\) moment equations. There are \(r^2\) multi-vertices and four mono-vertices. The rank follows from Theorem~5.2 and the rank-nullity theorem.
\end{proof}

\section{The lowest-degree example in the family}

For \(d=2\), we have \(r=4\) and
\[
\Omega=[0,7]\times[0,5].
\]
The cross-cuts are \(x=1\) and \(x=6\). The vertical \(\Tl\)-edges are
\[
V_j=\{j+1\}\times[1,4],
\qquad j=1,2,3,4,
\]
and the horizontal \(\Tl\)-edges are
\[
H_1=[1,6]\times\{1\},
\qquad
H_4=[1,6]\times\{4\},
\]
\[
H_2=[2,5]\times\{2\},
\qquad
H_3=[2,5]\times\{3\}.
\]
The relevant counts are
\[
n(H_1)=n(H_4)=6,
\qquad
n(H_2)=n(H_3)=n(V_j)=4.
\]
Thus there are eight \(\Tl\)-edges and twenty component vertices. The conformality matrix has size \(24\times20\) and rank \(19\).

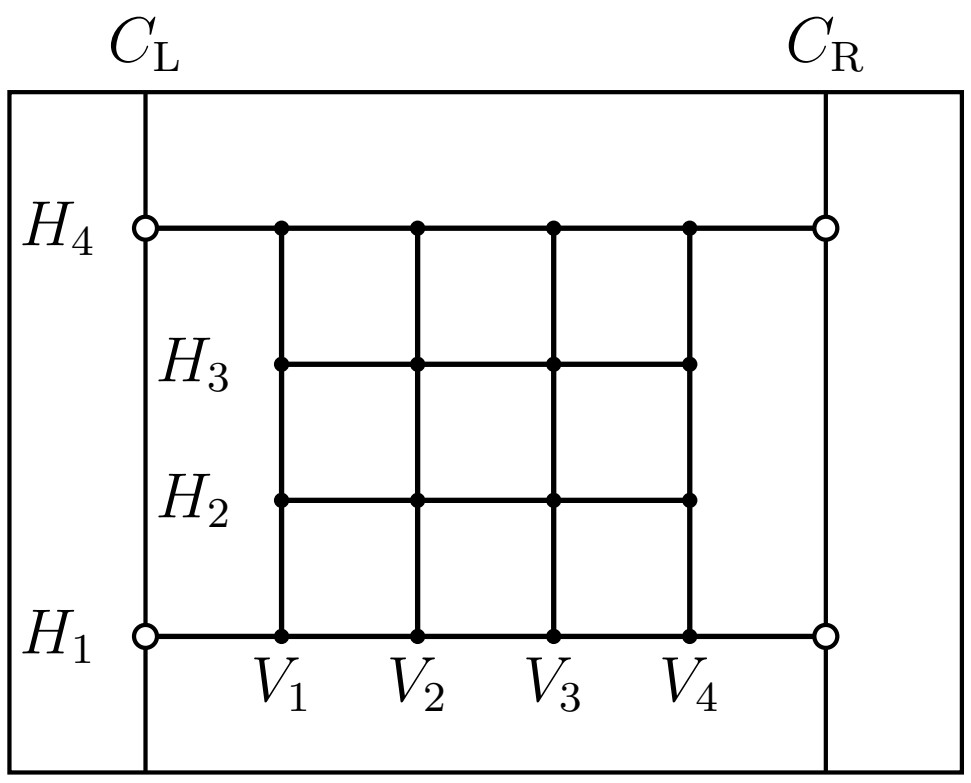
\begin{figure}[htbp]
\centering
\begin{tikzpicture}[scale=0.9]
  \draw[line width=0.8pt] (0,0) rectangle (7,5);

  \draw[line width=0.8pt] (1,0)--(1,5);
  \draw[line width=0.8pt] (6,0)--(6,5);

  \draw[line width=1pt] (1,1)--(6,1);
  \draw[line width=1pt] (1,4)--(6,4);

  \draw[line width=1pt] (2,2)--(5,2);
  \draw[line width=1pt] (2,3)--(5,3);

  \foreach \x in {2,3,4,5}{
    \draw[line width=1pt] (\x,1)--(\x,4);
  }

  \foreach \x in {2,3,4,5}{
    \foreach \y in {1,2,3,4}{
      \fill (\x,\y) circle (1.6pt);
    }
  }

  \foreach \x in {1,6}{
    \foreach \y in {1,4}{
      \draw[fill=white,line width=0.8pt]
        (\x,\y) circle (2.4pt);
    }
  }

  \node at (1,5.35) {$C_{\mathrm L}$};
  \node at (6,5.35) {$C_{\mathrm R}$};

  \node[anchor=east] at (0.8,1) {$H_1$};
  \node[anchor=east] at (1.8,2) {$H_2$};
  \node[anchor=east] at (1.8,3) {$H_3$};
  \node[anchor=east] at (0.8,4) {$H_4$};

  \foreach \x/\j in {2/1,3/2,4/3,5/4}{
    \node at (\x,0.65) {$V_{\j}$};
  }
\end{tikzpicture}
\caption{The \(T\)-mesh \(\mathcal T_2\). Filled circles are multi-vertices of its \(T\)-connected component. Open circles are its four mono-vertices. The two full-height interior lines are cross-cuts.}
\label{fig:T2}
\end{figure}

Since both coordinate sets are equally spaced with unit spacing, one may use
\[
u=(-1,3,-3,1)^T
\]
as a generator of each relevant Vandermonde null space. Accordingly, an intersection-coordinate generator is
\[
uu^T
=
\begin{pmatrix}
1&-3&3&-1\\
-3&9&-9&3\\
3&-9&9&-3\\
-1&3&-3&1
\end{pmatrix},
\]
extended by four zero mono-vertex coordinates.

This is the lowest-degree member of the present family. No claim is made that it minimizes the number of edges or vertices among all possible counterexamples.

\section{Discussion}

The example distinguishes constant rank from the existence of a diagonalizing ordering. Its conformality matrix has coordinate-independent rank throughout the specified structural similarity class, despite having unavoidable row dependencies.

The complete bipartite intersection pattern is central. It forces every vertical edge to meet the same set of horizontal supporting lines and every horizontal edge to meet the same set of vertical supporting lines. The vertical constraints first force the intersection matrix to have the form \(a\lambda^T\). One ordinary horizontal edge then forces \(\lambda\in\operatorname{span}\{b\}\). The remaining four mono-vertex variables vanish.

The failure of diagonalizability has a separate combinatorial cause. Every component edge has fewer than \(d+1\) mono-vertices. Thus no edge can occupy the last position of a diagonalizing ordering.

The restriction \(d\ge 2\) is necessary for this last-edge argument. For \(d=1\), two remaining mono-vertices meet the threshold \(d+1=2\), and the argument no longer applies.

Under Definitions~2.1--2.4, the implication
\[
\text{dimensionally absolutely stable}
\Longrightarrow
\text{diagonalizable}
\]
is therefore false. Any application to another notion of generalized component must first verify the common-coordinate requirement described in Remark~2.5.

\section*{Declaration of competing interest}

The authors declare that they have no known competing financial interests or personal relationships that could have appeared to influence the work reported in this paper.

\section*{Data availability}

No external datasets were used. All constructions and calculations needed to establish the results are included in the article.

\bibliographystyle{elsarticle-num}
\bibliography{references}

\end{document}